\documentclass[12pt]{article}
\font\smallit=cmti10 \font\smalltt=cmtt10 \font\smallrm=cmr9

\usepackage{amssymb,latexsym}

\makeatletter

\renewcommand\section{\@startsection {section}{1}{\z@}
{-30pt \@plus -1ex \@minus -.2ex} {2.3ex \@plus.2ex}
{\normalfont\normalsize\bfseries}}

\renewcommand\subsection{\@startsection{subsection}{2}{\z@}
{-3.25ex\@plus -1ex \@minus -.2ex} {1.5ex \@plus .2ex}
{\normalfont\normalsize\bfseries}}

\renewcommand{\@seccntformat}[1]{\csname the#1\endcsname. }
\makeatother
\newtheorem{defn}{Definition}
\newtheorem{prop}{Proposition}
\newtheorem{thm}{Theorem}
\newtheorem{cor}{Corollary}
\newcommand{\z}{\mathcal{Z}}
\newcommand{\n}{\mathcal{N}}

\begin{document}
\vspace*{-40pt} \centerline{\smalltt IJMC: \smallrm INTERNATIONAL JOURNAL OF MATHEMATICAL COMBINATORICS \smalltt(2010), Vol 2 \#A07-14}
\vskip 40pt

\begin{center}
\uppercase{\bf Generalizations of poly-Bernoulli numbers and
polynomials} \vskip 20pt
{\bf Hassan Jolany and M.R.Darafsheh, R.Eizadi Alikelaye}\\
{\smallit School of Mathematics,college of science,university of Tehran, Tehran, Iran}\\
{\tt Hassan.jolany@khayam.ut.ac.ir}\\ 
\end{center}
\vskip 30pt
\centerline{\smallit Received April 12, 2010: , Accepted: May 28, 2010 } 
\centerline{{\smallit Keywords}:Poly-Bernoulli polynomials,Euler
polynomials}
\centerline{,generalized Euler
polynomials,generalized poly-Bernoulli polynomials} \vskip 30pt

\centerline{\bf Abstract}

\noindent The Concepts of poly-Bernoulli numbers $B_n^{(k)}$ ,
poly-Bernoulli polynomials $B_n^{k}{(t)}$ and the generalized
poly-bernoulli numbers $B_{n}^{(k)}(a,b)$ are generalized to
$B_{n}^{(k)}(t,a,b,c)$ which is called the generalized
poly-Bernoulli polynomials depending on real parameters
\textit{a,b,c} .Some properties of these polynomials and some
relationships between $B_n^{k}$ , $B_n^{(k)}(t)$ ,
$B_{n}^{(k)}(a,b)$ and $B_{n}^{(k)}(t,a,b,c)$ are established.

\pagestyle{myheadings}

\thispagestyle{empty} \baselineskip=15pt \vskip 30pt

\section{Introduction}
In this paper we shall develop a number of generalizations of the
poly-Bernoulli numbers and polynomials , and obtain some results
about these generalizations.They are fundamental objects in the
theory of special functions.

Euler numbers are denoted with $B_k$ and are the coefficients of
Taylor expansion of the function $\frac{t}{e^{t}-1}$ ,as
following;
\[
\frac{t}{e^{t}-1}=\sum_{k=0}^{\infty}B_k\frac{t^k}{k!}. \]

The Euler polynomials $ E_n(x)$ are expressed in the following
series
\[
\frac{2e^{xt}}{e^t+1}=\sum_{k=0}^{\infty}E_k(x)\frac{t^k}{k!} .\]
for more details, see $[1]$-$[4].$

In $[10]$, Q.M.Luo, F.Oi and L.Debnath defined the generalization
of Euler polynomials\\
$E_k(x,a,b,c)$ which are expressed in the following series :

\[
\frac{2c^{xt}}{b^t+a^t}=\sum_{k=0}^{\infty}E_k(x,a,b,c)\frac{t^k}{k!}
.\]

where $a,b,c \in \z^+ $.They proved that\\

$I)$ for $a=1$ and $b=c=e$

\begin{equation}
E_k(x+1)=\sum_{j=0}^k\left({\begin{array}{l} k \\
j \end{array} }\right) E_j(x),
\end{equation}
and
\begin{equation}
E_k(x+1)+E_k(x)=2x^k,
\end{equation}

$II)$ for $a=1$ and $b=c$ ,
\begin{equation}
E_k(x+1,1,b,b)+E_k(x,1,b,b)=2x^k(\ln b)^k.
\end{equation}

In$[5]$, Kaneko introduced and studied poly-Bernoulli numbers
which generalize the classical Bernoulli numbers. Poly-Bernoulli
numbers $B_n^{(k)}$ with $k \in \z$ and $n \in \n$, appear in the
following power series:

\[
\frac{Li_k(1-e^{-x})}{1-e^{-x}}=\sum_{n=0}^{\infty}B_n^{(k)}\frac{t^n}{n!}
,\hspace{5 mm}(*)\]

where $k\in \z$ and

\[
Li_k(z)=\sum_{m=1}^{\infty}\frac{z^m}{m^k} .\hspace{6 mm}|z|<1\]

so for $k \leq 1$,

\[
Li_1(z)=-\ln(1-z),Li_0(z)=\frac{z}{1-z},Li_{-1}=\frac{z}{(1-z)^2},...
\hspace{2 mm}.\]

Moreover when $k \geq 1$, the left hand side of $(*)$ can be
wrriten in the form of "interated integrals"

\[
e^t\frac{1}{e^t-1}
\int_0^t\frac{1}{e^t-1}\int_0^t...\frac{1}{e^t-1}\int_0^t\frac{t}{e^t-1}dtdt...dt
\]

\[
=\sum_{n=0}^{\infty}B_n^{(k)}\frac{t^n}{n!} .\]

In the special case, one can see $B^{(1)}_n=B_n$.

\begin{defn}
The poly-Bernoulli polynomials,$B_n^{(k)}(t)$, are appeared in the
expansion of $\frac{Li_k(1-e^{-x})}{1-e^{-x}}e^{xt}$ as following,
\begin{equation}
\frac{Li_k(1-e^{-x})}{1-e^{-x}}e^{xt}=\sum_{n=0}^{\infty}\frac{B_n^{(k)}(t)}{n!}x^n
\end{equation}
for more details, see $[6]-[11]$.
\end{defn}

\begin{prop}
(Kaneko theorem$[6]$) The Poly-Bernoulli numbers of non-negative
index $k$ ,satisfy the following

\begin{equation}
B_n^{(k)}=(-1)^n\sum_{m=1}^{n+1}\frac{(-1)^{m-1}(m-1)!\left\{{\begin{array}{l}n
\\
m-1
\end{array}}\right\}}{m^k}
,\end{equation}

and for negative index k, we have

\begin{equation}
B_n^{(-k)}=\sum_{j=0}^{min(n,k)}(j!)^2\left\{{\begin{array}{l}n+1
\\
j+1\end{array}} \right\}\left\{{\begin{array}{l}k+1
\\
j+1 \end{array}}\right\} ,\end{equation}

where
\begin{equation}
\left\{{\begin{array}{l}n
\\
m
\end{array}
}\right\}
=\frac{(-1)^m}{m!}\sum_{l=0}^m(-1)^l\left(\begin{array}{l}m
\\
l
\end{array}\right)l^n \hspace{3 mm}m,n\geq0
\end{equation}
\end{prop}

\begin{defn}
Let $a,b>0$  and  $a\ne b$. The generalized poly-Bernoulli numbers
$B_n^{(k)}(a,b)$, the generalized poly-Bernoulli polynomials
$B_n^{(k)}(t,a,b)$ and the polynomial $B_n^{(k)}(t,a,b,c)$ are
appeared in the following series respectively;

\begin{equation}
\frac{Li_k(1-(ab)^{-t})}{b^t-a^{-t}}=\sum_{n=0}^{\infty}\frac{B_n^{(k)}(a,b)}{n!}t^n\hspace{4
mm}|t|< \frac{2\pi}{|\ln a+\ln b|} ,\end{equation}
\begin{equation}
\frac{Li_k(1-(ab)^{-t})}{b^t-a^{-t}}e^{xt}=\sum_{n=0}^{\infty}\frac{B_n^{(k)}(x,a,b)}{n!}t^n\hspace{4
mm} |t|<\frac{2\pi}{|\ln a+\ln b|},
\end{equation}

\begin{equation}
\frac{Li_k(1-(ab)^{-t})}{b^t-a^{-t}}c^{xt}=\sum_{n=0}^{\infty}\frac{B_n^{(k)}(x,a,b,c)}{n!}t^n\hspace{4
mm} |t|<\frac{2\pi}{|\ln a+\ln b|} ,\end{equation}
\end{defn}
\section{The main theorems}
We present some recurrence formulae for generalized poly-Bernoulli
polynomials.

\begin{thm}
Let $x \in \mathbb{R}$ and $n \geq 0$. For every positive real
numbers a,b and c such that $a\ne b$ and $b>a$ , we have\\
\begin{equation}
B_n^{(k)}(a,b)=B_n^{(k)}\left(\frac{-\ln b}{\ln a+\ln
b}\right)(\ln a+\ln b)^n,
\end{equation}
\begin{equation}
B_j^{(k)}(a,b)=\sum_{i=1}^{j}(-1)^{j-i}(\ln a+\ln b)^i(\ln
b)^{j-i}\left(\begin{array}{l}j
\\
i
\end{array}\right)B_j^{(k)}
,\end{equation}
\begin{equation}
B_n^{(k)}(x;a,b,c)=\sum_{l=0}^n\left(\begin{array}{l} n
\\
l
\end{array}\right)(\ln
c)^{n-l}B_l^{(k)}(a,b)x^{n-l} ,\end{equation}

\begin{equation}
B_n^{(k)}(x+1;a,b,c)=B_n^{(k)}(x;ac,\frac{b}{c},c) ,\end{equation}

\begin{equation}
B_n^{(k)}(t)=B_n^{(k)}(e^{t+1},e^{-t}) ,\end{equation}
\begin{equation}
B_n^{(k)}(x,a,b,c)=(\ln a+\ln b)^nB_n^{(k)}(\frac{-\ln b+x\ln
c}{\ln a+\ln b}) .\end{equation}
\end{thm}

{\bf Proof:}By applying definition 2, we have

(11)    \[
\frac{Li_k(1-(ab)^{-t})}{b^t-a^{-t}}=\sum_{n=0}^{\infty}\frac{B_n^{(k)}(a,b)}{n!}t^n
\]

\[=
\frac{1}{b^t}\left(\frac{Li_k(1-e^{-t\ln ab})}{1-e^{-t\ln
ab}}\right)=e^{-t\ln b}\left(\frac{Li_k(1-e^{-t\ln
ab})}{1-e^{-t(\ln ab)}}\right)
\]

\[
=\sum_{n=0}^{\infty}B_n^{(k)}\left(\frac{-\ln b}{\ln a+\ln
b}\right)(\ln a+\ln b)^n\frac{t^n}{n!}
\]

Therefore
\[
B^{(k)}_n(a,b)=B_n^{(k)}\left(\frac{-\ln b}{\ln a+\ln
b}\right)(\ln a+\ln b)^n .\] (12)
\[
\frac{Li_k(1-(ab)^{-t})}{b^t-a^{-t}}=\frac{1}{b^t}\left(\frac{Li_k(1-(ab)^{-t\
})}{1-e^{-t\ln ab}}\right) =\left(\sum_{k=0}^{\infty}\frac{(\ln
b)^k}{k!}(-1)^kt^k\right)\left(\sum_{n=0}^{\infty}B_n^{(k)}\frac{(\ln
a+\ln b)^n}{n!}t^n\right)\]

\[=
\sum_{j=0}^{\infty}\left(\sum_{i=0}^j(-1)^{j-i}B_i^{(k)}\frac{(\ln
a+\ln b)^i}{i!(j-i)!}(\ln b)^{j-i}\right)t^j
\]

So we have
\[B_j^{(k)}(a,b)=\sum_{i=0}^{j}(-1)^{j-i}(\ln a+\ln b)^i(\ln b)^{j-i}\left(\begin{array}{l}j
\\
i
\end{array}\right)B_i^{(k)}
.\]

(13)
\[\frac{Li_k(1-(ab)^{-t})}{b^t-a^{-t}}c^{xt}=
\sum_{n=0}^{\infty}B_n^{(k)}(x,a,b,c)\frac{t^n}{n!}
\]

\[=
\left(\sum_{l=0}^{\infty}B_l^{(k)}(a,b)
\frac{t^l}{l!}\right)\left(\sum_{i=0}^{\infty}\frac{(\ln
c)^it^i}{i!}x^i\right)
\]

\[=
\sum_{l=0}^{\infty}\sum_{i=0}^{l}\frac{(\ln c
)^{l-i}}{i!(l-i)!}B_i^{(k)}(a,b)x^{l-i}t^l
\]

\[=
\sum_{n=0}^{\infty}\left(\sum_{l=0}^n\left(\begin{array}{l}n
\\
l
\end{array}\right)(\ln
c)^{n-l}B_l^{(k)}(a,b)x^{n-l}\right)\frac{t^n}{n!} .\]

(14)

\[
\frac{Li_k(1-(ab)^{-t})}{b^t-a^{-t}}c^{(x+1)t}=
\frac{Li_k(1-(ab)^{-t})}{b^t-a^{-t}}c^{xt}.c^t
\]

\[=
\frac{Li_k(1-(ab)^{-t})}{\left(\frac{b}{c}\right)^t-(ac)^{-t}}c^{xt}=
\sum_{n=0}^{\infty}B_n^{(k)}(x;ac,\frac{b}{c},c)\frac{t^n}{n!} .\]
(15)
\[
\frac{Li_k(1-e^{-x})}{1-e^{-x}}e^{xt}=
\frac{Li_k(1-e^{-x})}{e^{-xt}-e^{-x-xt}}=
\frac{Li_k(1-e^{-x})}{(e^{-t})^x-(e^{1+t})^{-x}}
\]

so,
\[
B_n^{(k)}(t)=B_n^{(k)}(e^{t+1},e^{-t}).
\]
(16)we can write
\[
\sum_{n=0}^{\infty}B_n^{(k)}(x,a,b,c)\frac{t^n}{n!}=
\frac{Li_k(1-(ab)^{-t})}{b^t-a^{-t}}c^{xt}=
\frac{1}{b^t}\frac{Li_k(1-(ab)^{-t})}{(1-(ab)^{-t})}c^{xt}
\]

\[=e^{t({-\ln b+x\ln c})}\left(\frac{Li_k(1-e^{-t\ln ab})}{1-e^{-t(\ln
ab)}}\right)= \sum_{n=0}^{\infty}(\ln a+\ln
b)^nB_n^{(k)}\left(\frac{-\ln b+x\ln c}{\ln a+\ln
b}\right)\frac{t^n}{n!}
\]

so
\[
B_n^{(k)}(x,a,b,c)=(\ln a+\ln b)^nB_n^{(k)}\left(\frac{-\ln b+x\ln
c}{\ln a+\ln b}\right) .\]
\\

\begin{thm}
Let $x \in \mathbb{R}$ , $n \geq 0$. For every positive real
numbers a,b such that  $a\ne b$ and $b>a>0$, we have

\begin{equation}
B_n^{(k)}(x+y,a,b,c)=\sum_{l=0}^{\infty}
\left(\begin{array}{l}n
\\
l
\end{array}\right)(\ln
c)^{n-l}B_l^{(k)}(x;a,b,c)y^{n-l}
\newline\\=\sum_{l=0}^n\left(\begin{array}{l}
n
\\
l
\end{array}\right)(\ln c)^{n -l}B_l^{(k)}(y,a,b,c)x^{n-l}
.\end{equation}

\end{thm}

{\bf Proof:} We have \[
\frac{Li_k(1-(ab)^{-t})}{b^t-a^{-t}}c^{(x+y)t}=\sum_{n=0}^{\infty}B_n^{(k)}(x+y;a,b,c)\frac{t^n}{n!}
\]

\[=\frac{Li_k(1-(ab)^{-t})}{b^t-a^{-t}}c^{xt}.c^{yt}=\left(\sum_{n=0}^{\infty}B_n^{(k)}(x;a,b,c)\frac{t^n}{n!}\right)
\left(\sum_{i=0}^{\infty}\frac{y^i(\ln c)^i}{i!}t^i\right)
\]
\[
=\sum_{n=0}^{\infty}\left(\sum_{l=0}^n\left(\begin{array}{l}n
\\
l
\end{array}\right)y^{n-l}(\ln
c)^{n-l}B_l^{(k)}(x,a,b,c)\right)\frac{t^n}{n!}
\]
so we get
\[
\frac{Li_k(1-(ab)^t)}{b^t-a^{-t}}c^{(x+y)t}=\frac{Li_k(1-(ab)^{-t}}{b^t-a^{-t}}c^{yt}c^{xt}
\]
\[
=\sum_{n=0}^{\infty}\left(\sum_{l=0}^n\left(\begin{array}{l}n
\\
l
\end{array}\right)x^{n-l}(\ln
c)^{n-l}B_l^{(k)}(y,a,b,c)\right)\frac{t^n}{n!} .\]
\\

\begin{thm}
 Let $x \in \mathbb{R}$ and $n \geq 0$.For every positive real
numbers a,b and c such that  $a\ne b$ and $b>a>0$, we have

\begin{equation}
B_n^{(k)}(x;a,b,c)=\sum_{l=0}^n\left(\begin{array}{l}n
\\
l
\end{array}\right)(\ln c)^{n-l}B_l^{(k)}\left(\frac{-\ln b}{\ln a+\ln
b}\right)(\ln a+\ln b)^lx^{n-l},
\end{equation}
\begin{equation}
B_n^{(k)}(x;a,b,c)=\sum_{l=0}^n\sum_{j=0}^l(-1)^{l-j}\left(\begin{array}{l}n
\\
l
\end{array}\right)\left(\begin{array}{l}
l
\\
j
\end{array}\right)(\ln c)^{n-l}(\ln b)^{l-j}(\ln a+\ln
b)^jB_j^{(k)}x^{n-k}.
\end{equation}
\end{thm}

{\bf Proof:}
\\$\hspace{10 mm}(18)$By applying Theorems 1 and 2 we know ,

\[
B_n^{(k)}(x;a,b,c)=\sum_{l=0}^n\left(\begin{array}{l}n
\\
l\end{array}\right)(\ln c)^{n-l}B_l^{(k)}(a,b)x^{n-l}
\]
and
\[
B_n^{(k)}(a,b)=B_n^{(k)}(\frac{-ln b}{\ln a+\ln b})(\ln a+ \ln
b)^n
\]

The relation (18)will follow if we combine these formulae.

(19) proof is similar to(18).
\\

Now,we give some results about derivatives and integrals of the
generalized poly-Bernoulli polynomials in the following theorem.

\begin{thm}
Let $x \in \mathbb{R}$.If $a,b$ and $c>0$  , $a \ne b$ and $b>a>0$
,For any non-negative integer l and real numbers $\alpha$ and
$\beta$ we have
\begin{equation}
\frac{\partial ^l B_n^{(k)}(x,a,b,c)}{\partial
x^l}=\frac{n!}{(n-l)!}(\ln c)^lB_{n-l}^{(k)}(x,a,b,c)
\end{equation}

\begin{equation}
\int_{\alpha}^{\beta}B_n^{(k)}(x,a,b,c)dx=\frac{1}{(n+1)\ln
c}[B_{n+1}^{(k)}(\beta,a,b,c)-B_{n+1}^{(k)}(\alpha,a,b,c)]
\end{equation}
\end{thm}

{\bf Proof:} By applying (20) and (21) can be proved by induction
on \textit{n}.

In [9], GI-Sang Cheon investigated the classical relationship
between Bernoulli and Euler polynomials, in this paper we study
the relationship between the generalized poly-Bernoulli and Euler
polynomials.

\begin{thm}
For $b>0$ we have
\[
B_n^{(k_1)}(x+y,1,b,b)=\frac{1}{2}\sum_{k=0}^n\left(\begin{array}{l}n
\\
k\end{array}\right)[B_n^{(k_1)}(y,1,b,b)+B_n^{(k_1)}(y+1,1,b,b)]E_{n-k}(x,1,b,b)
.\]
\end{thm}

{\textbf Proof:}We know
\[
B_n^{(k_1)}(x+y,1,b,b)=\sum_{k=0}^{\infty}\left(\begin{array}{l}n
\\
k\end{array}\right)(\ln b)^{n-k}B_k^{(k_1)}(y;1,b,b)x^{n-k}
\]

and
\[
E_k(x+y,1,b,b)+E_k(x,1,b,b)=2x^k(\ln b)^k
\]

So, We obtain $B_n^{(k_1)}(x+y,1,b,b)=$

\[
\frac{1}{2}\sum_{k=0}^n\left(\begin{array}{l} n
\\
k
\end{array}\right)(\ln b)^{n-k}B_k^{(k_1)}(y;1,b,b)
\left[\frac{1}{(\ln
b)^{n-k}}(E_{n-k}(x,1,b,b)+E_{n-k}(x+1,1,b,b))\right ]\
\]

\[
=\frac{1}{2}\sum_{k=0}^n\left(\begin{array}{l}n
\\
k
\end{array}\right)B_k^{(k_1)}(y;1,b,b)\left[\left(E_{n-k}(x,1,b,b)+\sum_{j=0}^{n-k}\left(\begin{array}{l}n-k
\\
j
\end{array}\right)E_j(x,1,b,b)\right)\right ]\
\]

\[
=\frac{1}{2}\sum_{k=0}^n\left(\begin{array}{l}n
\\
k\end{array}\right)B_k^{(k_1)}(y;1,b,b)E_{n-k}(x,1,b,b)
\]
\[
+\frac{1}{2}\sum_{j=0}^n\left(\begin{array}{l} n
\\
j
\end{array}\right)E_j(x;1,b,b)\sum_{k=0}^{n-j}\left(
\begin{array}{l}
n-j
\\
k
\end{array}\right)B_k^{(k_1)}(y,1,b,b)
\]
\[
=\frac{1}{2}\sum_{k=0}^n\left(\begin{array}{l}n
\\
k\end{array}\right)B_k^{(k_1)}(y;1,b,b)E_{n-k}(x,1,b,b)+
\]
\[
\frac{1}{2}\sum_{j=0}^n\left(\begin{array}{l} n
\\
j
\end{array}\right)B_{n-j}^{(k_1)}(y+1;1,b,b)E_j(x,1,b,b)
\]

So we have
\[
B_n^{(k_1)}(x+y,1,b,b)=\frac{1}{2}\sum_{k=0}^n\left(\begin{array}{l}n
\\
k\end{array}\right)[B_n^{(k_1)}(y,1,b,b)+B_n^{(k_1)}(y+1,1,b,b)]E_{n-k}(x,1,b,b)
.\]
\begin{cor}
In theorem 5, if we set $k_1=1$ and $b=e$, we obtain
\[
B_n(x)=\sum_{(k=0),(k\ne1)}^n\left(\begin{array}{l}n
\\
k\end{array}\right)B_k E_{n-k}(x) .\]
\end{cor}
For more detail see [7].

{\bf Acknowledgement:}The authors would like to thank the referee
for has valuable comments concerning theorem 5.

\bibliographystyle{plain}

\end{document}